\documentclass[12pt]{amsart}

\newtheorem{theorem}{Theorem}[section]
\newtheorem{proposition}[theorem]{Proposition}
\newtheorem{lemma}[theorem]{Lemma}
\newtheorem{corollary}[theorem]{Corollary}

\numberwithin{equation}{section}

\begin{document}

\title[K\"ahler structures over symmetric products]{On K\"ahler
structures over symmetric products of a Riemann surface}

\author{Indranil Biswas}

\address{School of Mathematics, Tata Institute of Fundamental
Research, Homi Bhabha Road, Bombay 400005, India}

\email{indranil@math.tifr.res.in}



\date{}

\begin{abstract}
Given a positive integer $n$, and a compact connected Riemann 
surface $X$,
we prove that the symmetric product $S^n(X)$ admits a K\"ahler form
of nonnegative holomorphic bisectional curvature if and only if
$\text{genus}(X)\, \leq\, 1$. If $n$ is greater than or equal to the 
gonality of $X$, we prove that $S^n(X)$ does not admit any K\"ahler
form of nonpositive holomorphic sectional curvature. In particular,
if $X$ is hyperelliptic, then $S^n(X)$ admits a K\"ahler
form of nonpositive holomorphic sectional curvature if and only if
$n\,=\,1\, \leq\, \text{genus}(X)$.
\end{abstract}

\maketitle

\section{Introduction}\label{sec1}

This work was inspired by \cite{BR}. The following is the main
theorem proved there (see \cite[Theorem 1.1]{BR}):

\textit{Let $X$ be a compact connected Riemann surface of genus
at least two. If
$$
n\, \leq\, 2({\rm genus}(X)-1)\, ,
$$
then the symmetric product $S^n(X)$ does not admit any
K\"ahler form of nonnegative holomorphic bisectional curvature.}

By nonnegative holomorphic bisectional curvature we mean that
all holomorphic bisectional curvatures are nonnegative. Similarly,
by nonpositive holomorphic sectional curvature we mean that 
all holomorphic sectional curvatures are nonpositive.

Our aim here is to settle the cases not considered in \cite{BR}.
We prove the following theorem (see Theorem \ref{thm1}):

\begin{theorem}\label{thm0}
Let $X$ be a compact connected Riemann surface, and let $n$ be a
fixed positive integer. The symmetric product $S^n(X)$ admits a
K\"ahler form of nonnegative holomorphic bisectional curvature
if and only if ${\rm genus}(X)\, \leq\, 1$.
\end{theorem}

The gonality of a compact connected Riemann surface $X$ is the
smallest integer $d$ such that there is a nonconstant holomorphic
map $X\, \longrightarrow\, {\mathbb C}{\mathbb P}^1$ of degree $d$.

We prove the following (see Proposition \ref{prop3}):

\begin{proposition}\label{prop0}
Let $X$ be a compact connected Riemann surface. Let $d$ be the
gonality of $X$. Take any integer $n\, \geq\, d$. Then
$S^n(X)$ does not admit any K\"ahler form
of nonpositive holomorphic sectional curvature.
\end{proposition}

Note that a K\"ahler form of nonpositive holomorphic bisectional 
curvature is also a K\"ahler form of nonpositive holomorphic 
sectional curvature.

Proposition \ref{prop0} gives the following (see
Corollary \ref{cor1}):

\begin{corollary}
Let $X$ be a hyperelliptic Riemann surface. The compact
complex manifold $S^n(X)$ admits a K\"ahler form of
nonpositive holomorphic sectional curvature if
and only if $n\, =\, 1\, \leq\, {\rm genus}(X)$.
\end{corollary}

\section{Nonnegative holomorphic bisectional curvature}

Let $X$ be a compact connected Riemann surface. The genus of
$X$ will be denoted by $g$. For any positive integer $n$, let
$S^n(X)$ be the symmetric product of $X$. Therefore, $S^n(X)$
is the quotient of the Cartesian product
$X^n$ by the natural action of the group of
permutations of the index set $\{1\, , \cdots \, ,n\}$. It
is known that $S^n(X)$ is a smooth complex projective variety of
dimension $n$.

In this section we address the question whether $S^n(X)$ admits
a K\"ahler structure of nonnegative holomorphic bisectional curvature.

\begin{proposition}\label{prop1}
Assume that $g\, \geq\,2$. 
The complex manifold $S^n(X)$ does not admit any K\"ahler
metric of nonnegative holomorphic bisectional curvature.
\end{proposition}

\begin{proof}
If $n \, \leq\, 2g-2$, then this is proved in
\cite{BR} (see \cite[Theorem 1.1]{BR}). So we assume that
$n \, > \, 2g-2$.

Let $J^n(X)$ be the Picard variety parametrizing isomorphism
classes of holomorphic line bundles over $X$ of degree $n$.
Fix a Poincar\'e line bundle
\begin{equation}\label{e1}
{\mathcal L} \, \longrightarrow\, X\times J^n(X)
\end{equation}
(see \cite[Ch.~IV, \S~2]{ACGH} for the construction of a
Poincar\'e line bundle). Let
\begin{equation}\label{e2}
p\, :\, X\times J^n(X)\, \longrightarrow\, J^n(X)
\end{equation}
be the projection to the second factor. Since
$n \, > \, 2g-2$, we have
$\text{degree}(L^*\otimes K_X)\, <\, 0$
for any line bundle $L$ of degree $n$ on $X$, where
$K_X$ is the holomorphic cotangent bundle of $X$. Hence
by Serre duality, we have
$$
H^1(X,\, L)\,=\,
H^0(X,\, L^*\otimes K_X)^*\,=\,0
$$
for any holomorphic line bundle $L$ of degree $n$.
Therefore, for the projection $p$ in \eqref{e2}, we have
$$
R^1p_* {\mathcal L} \, =\, 0\, ,
$$
and also
\begin{equation}\label{e3}
{\mathcal V} \, :=\, p_* {\mathcal L}\,\longrightarrow\, J^n(X)
\end{equation}
is a holomorphic vector bundle of rank $n-g+1$.

Consider the projective bundle
\begin{equation}\label{e3a}
f\, :\, P({\mathcal V})\, \longrightarrow\, J^n(X)
\end{equation}
parametrizing the lines in the fibers of the vector bundle
${\mathcal V}$ in \eqref{e3}. Since any two Poincar\'e line bundles
over $X\times J^n(X)$ differ by tensoring with a line bundle
pulled back from $J^n(X)$ \cite[p. 166]{ACGH}, using
the projection formula we conclude that the projective bundle
$P({\mathcal V})$ is actually independent of the choice of the
Poincar\'e line bundle $\mathcal L$.

The total space of $P({\mathcal V})$ is identified with
$S^n(X)$. Points of $P({\mathcal V})$ parametrize
isomorphism classes of pairs
of the form $(L\, ,s)$, where $L\, \in\, J^n(X)$ and
$s$ is a holomorphic section of $L$ which is not identically
zero. The identification of $P({\mathcal V})$ with $S^n(X)$ sends
any pair $(L\, ,s)$ to the divisor of the section $s$.

Assume that $S^n(X)$ has a K\"ahler
metric of nonnegative holomorphic bisectional curvature. Therefore,
the total space of $P({\mathcal V})$ has a K\"ahler
metric of nonnegative holomorphic bisectional curvature. This implies
that the anticanonical line bundle $K^{-1}_{P({\mathcal V})}\, :=\,
\bigwedge^n TP({\mathcal V})$ is numerically effective.

Let $T_{\rm rel}\, \subset\, TP({\mathcal V})$ be the relative
tangent bundle for the projection $f$ in \eqref{e3a}. In other words,
$T_{\rm rel}$ is the kernel of the differential $df\, :\,
TP({\mathcal V})\, \longrightarrow\, f^* TJ^n(X)$. Since
the tangent bundle $TJ^n(X)$ is trivial, the holomorphic
line bundle
\begin{equation}\label{dl}
\det (T_{\rm rel})\, :=\, \bigwedge\nolimits^{n-g} T_{\rm rel}
\, \longrightarrow\, P({\mathcal V})
\end{equation}
is identified with $K^{-1}_{P({\mathcal V})}$. It was noted above
that $K^{-1}_{P({\mathcal V})}$ is numerically effective. Therefore,
we now conclude that the line bundle $\det (T_{\rm rel})$
in \eqref{dl} is numerically effective.

Since $\det (T_{\rm rel})$ is numerically effective, the vector
bundle $\mathcal V$ has the property that
\begin{equation}\label{e4}
c_2({\rm ad}({\mathcal V})) \, =\, 0
\end{equation}
\cite[Theorem 1.1]{BB}, where ${\rm ad}({\mathcal V})\, \subset\, 
End({\mathcal V})\,=\, {\mathcal V}\otimes{\mathcal V}^*$ is the 
subbundle of co-rank one defined by the sheaf of endomorphisms of 
trace zero.

For a particular choice of the Poincar\'e bundle $\mathcal L$,
\begin{equation}\label{e5}
ch({\mathcal V}) \, =\, (n-g+1) - \theta\, ,
\end{equation}
where $\theta\, \in\, H^2(J^n(X), \, {\mathbb Q})$
is the class of a theta divisor; see lines 2--4 (from top)
of \cite[p. 336]{ACGH}. From \eqref{e5} it follows immediately
that
$$
c_2({\rm ad}({\mathcal V})) \, =\,
- ch_2({\rm ad}({\mathcal V}))\, =\,\theta^2\, .
$$
Note that $\theta^2\, \not=\, 0$ because $\theta$ is an
ample class and $\dim J^n(X)\, \geq\, 2$. But this contradicts
\eqref{e4}. Therefore, we conclude that $S^n(X)$ does not admit any 
K\"ahler metric of nonnegative holomorphic bisectional curvature.
\end{proof}

\begin{proposition}\label{prop2}
Assume that $g\, \leq\, 1$.
The complex manifold $S^n(X)$ admits a K\"ahler
metric of nonnegative holomorphic bisectional curvature.
\end{proposition}

\begin{proof}
First assume that $g\, =\, 0$. Then $S^n(X)$ is biholomorphic
to ${\mathbb C}{\mathbb P}^n$. A Fubini--Study metric on
${\mathbb C}{\mathbb P}^n$ has positive holomorphic bisectional
curvature. (In fact this property of existence of a
K\"ahler form of positive holomorphic bisectional
curvature characterizes
${\mathbb C}{\mathbb P}^n$ \cite{Mo}, \cite{SY}.)

Now assume that $g\, =\, 1$. Consider the vector
bundle ${\mathcal V}\, \longrightarrow\, J^n(X)$
constructed as in \eqref{e3}. As before, we have
$$
P({\mathcal V})\, =\, S^n(X)\, .
$$

We will show that the projective bundle
$P({\mathcal V})$ over $J^n(X)$ is given by a representation of
$\pi_1(J^n(X), x_0)$ in the projective unitary group $\text{PU}(n)$.

To prove this, first note that a theorem due to Ein--Lazarsfeld and
Kempf says that the vector bundle $\mathcal V$ on $J^n(X)$ is stable 
\cite[p. 149, Theorem]{EL},
\cite[p. 285, Theorem 3]{Ke}. Since $\dim J^n(X)\, =\, 1$,
any stable vector bundle on $J^n(X)$ has a projectively flat
unitary connection \cite{NS}; this also follows from
Atiyah's classification of vector bundles over an elliptic
curve \cite{At}. This proves the above statement that
$P({\mathcal V})$ is given by a representation of
$\pi_1(J^n(X), x_0)$ in $\text{PU}(n)$. In fact, Kempf
showed that the natural $L^2$--metric on $\mathcal V$ is
projectively flat \cite[p. 286, Corollary 6]{Ke}.

Fix a representation
$$
\rho\, :\, \pi_1(J^n(X), x_0)\, \longrightarrow\, \text{PU}(n)
$$
corresponding to a projectively flat unitary connection on $P({\mathcal 
V})$. Let $\widetilde{J^n(X)}\,\longrightarrow\, J^n(X)$ be the
universal cover associated to the base point $x_0$; note that
$\widetilde{J^n(X)}$ is a principal $\pi_1(J^n(X), x_0)$--bundle
over $J^n(X)$. Let
\begin{equation}\label{vp}
\varphi\, :\, F_\rho\, :=\, \widetilde{J^n(X)}\times^\rho 
{\mathbb C}{\mathbb P}^{n-1}\, \longrightarrow\, J^n(X)
\end{equation}
be the holomorphic fiber bundle associated to
the principal $\pi_1(J^n(X), x_0)$--bundle
$$\widetilde{J^n(X)}\,\longrightarrow\, J^n(X)$$ for the action of
$\pi_1(J^n(X), x_0)$ on ${\mathbb C}{\mathbb P}^{n-1}$ defined by
$\rho$ and the standard action of $\text{PU}(n)$ on ${\mathbb 
C}{\mathbb P}^{n-1}$. The holomorphic fiber bundle
$F_\rho\, \stackrel{\varphi}{\longrightarrow}\, J^n(X)$ is equipped
with a flat holomorphic connection; this connection will be denoted
by $\nabla$. The total space of $P({\mathcal V})$
is biholomorphic to the total space of $F_\rho$ because both the
fiber bundles are given by $\rho$.

Fix a Fubini--Study metric $\omega_{FS}$
on ${\mathbb C}{\mathbb P}^{n-1}$ preserved by the standard
action of $\text{PU}(n)$ on ${\mathbb C}{\mathbb P}^{n-1}$.
This metric defines a closed $(1\, ,1)$--form $\omega_F$ on the total
space of $F_\rho$ which is
K\"ahler on the fibers of the projection $\varphi$
in \eqref{vp}; this form
$\omega_F$ is uniquely determined by the following conditions:
\begin{enumerate}
\item for any point $z\, \in\, F_\rho$ and any horizontal tangent vector
$v\, \in\,T^{1,0}_z F_\rho\oplus T^{0,1}_z F_\rho$ for the above connection
$\nabla$,
$$
i_v \omega_F(z) \, =\, 0
$$
(it is the contraction of $\omega_F(z)\, \in\, \Omega^{1,1}_z
F_\rho$ by $v$), and

\item the restriction of $\omega_F$ to each fiber of
the projection $F_\rho\, \longrightarrow\, J^n(X)$ is $\omega_{FS}$.
\end{enumerate}

The closed $(1\, ,1)$--form $\omega_F$ can also be described as follows.
Consider the $(1\, ,1)$--form $p^*_{{\mathbb C}{\mathbb P}^{n-1}}\omega_{FS}$
on $\widetilde{J^n(X)}\times{\mathbb C}{\mathbb P}^{n-1}$, where
$p_{{\mathbb C}{\mathbb P}^{n-1}}$ is the projection of
$\widetilde{J^n(X)}\times{\mathbb C}{\mathbb P}^{n-1}$ to ${\mathbb C}{\mathbb 
P}^{n-1}$. This form $p^*_{{\mathbb C}{\mathbb P}^{n-1}}\omega_{FS}$ descends
to the quotient space $\widetilde{J^n(X)}\times^\rho{\mathbb C}{\mathbb 
P}^{n-1}$. This descended form coincides with $\omega_F$.

For any K\"ahler form $\omega_0$ on $J^n(X)$, it is 
straight--forward to check that
$$
\omega_F+ \varphi^*\omega_0
$$
is a K\"ahler form on the total space of $F_\rho$, where 
$\varphi$ is the projection in \eqref{vp}.

Since $J^n(X)$ is an elliptic curve, it has a K\"ahler metric
$\omega_{J^n(X)}$ of curvature zero (take any translation 
invariant metric on $J^n(X)$).

Since the Fubini--Study metric $\omega_{FS}$ has positive holomorphic
bisectional curvature, and the curvature of $\omega_{J^n(X)}$ is
zero, we conclude that the holomorphic
bisectional curvature of the K\"ahler form
$$
\omega_F+ \varphi^*\omega_{J^n(X)}
$$
on the total space of $F_\rho$
is nonnegative. Indeed, this follows immediately from the fact 
that the K\"ahler form $\omega_F+\varphi^*\omega_{J^n(X)}$ is 
locally the 
product of $\omega_{FS}$ on ${\mathbb C}{\mathbb P}^{n-1}$ and
$\omega_{J^n(X)}$ on $J^n(X)$. Since the
total space of $F_\rho$ is biholomorphic to the total space of
$P({\mathcal V})$, the proof of the proposition is complete.
\end{proof}

Combining Proposition \ref{prop1} and Proposition \ref{prop2},
we have the following theorem:

\begin{theorem}\label{thm1}
Let $X$ be a compact connected Riemann surface, and let $n$ be
a fixed positive integer. The compact complex manifold $S^n(X)$ 
admits a K\"ahler
form of nonnegative holomorphic bisectional curvature if and
only if ${\rm genus}(X)\, \leq\, 1$.
\end{theorem}

\section{Nonpositive holomorphic sectional curvature}

We begin by recalling a property of compact K\"ahler
manifolds with nonpositive holomorphic sectional curvature.

\begin{lemma}\label{lem1}
Let $M$ be a compact connected K\"ahler manifold admitting a
K\"ahler form of nonpositive holomorphic sectional curvature.
Then there is no nonconstant holomorphic map from ${\mathbb C}
{\mathbb P}^1$ to $M$.
\end{lemma}

\begin{proof}
See \cite[p. 40, Corollary 4.5]{Gr} for a proof of this lemma.
In \cite{Gr}, ``negatively curved'' K\"ahler metric means one
with nonpositive holomorphic sectional curvature (see Definition
in \cite[p. 39]{Gr}).
\end{proof}

Let $X$ be a compact connected Riemann surface. 
The \textit{gonality} of $X$ is the
smallest integer $d$ such that there is a nonconstant holomorphic
map $X\, \longrightarrow\, {\mathbb C}{\mathbb P}^1$ of degree $d$.
Equivalently, gonality of $X$
is the smallest integer $d$ such that there
is a holomorphic line bundle $\xi$ on $X$ of degree $d$ with
$\dim H^0(X,\, \xi)\, \geq\, 2$. Therefore, from the Riemann--Roch
theorem it follows that the gonality of $X$ is bounded above by
$\text{genus}(X)+1$.

If the gonality of $X$ is two, then
$X$ is called a hyperelliptic Riemann surface.

\begin{proposition}\label{prop3}
Let $X$ be a compact connected Riemann surface of gonality $d$.
If $n\, \geq\, d$, then $S^n(X)$ does not admit any K\"ahler form
of nonpositive holomorphic sectional curvature. In particular,
$S^n(X)$ does not admit a K\"ahler form
of nonpositive holomorphic sectional curvature if
$n\, > {\rm genus}(X)$.
\end{proposition}

\begin{proof}
Let $\phi\, :\, X \, \longrightarrow\, {\mathbb C}{\mathbb P}^1$
be a nonconstant holomorphic map of degree $\delta$. Then
we get a holomorphic map
$$
\widetilde{\phi}\, :\, {\mathbb C}{\mathbb P}^1
 \, \longrightarrow\, S^\delta(X)
$$
that sends any $x\, \in\, {\mathbb C}{\mathbb P}^1$ to the
scheme--theoretic inverse image $\phi^{-1}(x)$ of $x$.
This map $\widetilde{\phi}$ is clearly nonconstant.

If
$$
\psi\, :\, {\mathbb C}{\mathbb P}^1
 \, \longrightarrow\, S^a(X)
$$
is a nonconstant holomorphic map, and $b\, >\, a$, then we
can construct a nonconstant holomorphic map
$$
{\mathbb C}{\mathbb P}^1\, \longrightarrow\, S^b(X)
$$
as follows: fix a point $y_0\,\in\, X$, and send any
$x\, \in\, {\mathbb C}{\mathbb P}^1$ to $\psi(x) +(b-a)y_0\, \in\,
S^b(X)$.

In view of these two observations, the proposition follows from
Lemma \ref{lem1}.
\end{proof}

\begin{corollary}\label{cor1}
Let $X$ be a hyperelliptic Riemann surface. The compact complex
manifold $S^n(X)$ admits a K\"ahler form of nonpositive 
holomorphic sectional curvature if and only if $n\, =\, 1$.
\end{corollary}

\begin{proof}
Since $X$ is hyperelliptic, we have ${\rm genus}(X)\, > \,0$,
hence $X\,=\, S^1(X)$ admits a K\"ahler form of nonpositive 
holomorphic sectional curvature. If $n\, \geq\, 2$, from
Proposition \ref{prop3} we know that
$S^n(X)$ does not admit any K\"ahler form of nonpositive
holomorphic sectional curvature,
\end{proof}


\end{document}